\documentclass[12pt,letterpaper]{amsart}    

\usepackage[parfill]{parskip}    

\usepackage{graphicx}

\usepackage{amssymb} 

\usepackage{amsmath} 

\usepackage{amsthm} 

\usepackage{mathrsfs}

\usepackage{epstopdf}

\usepackage{enumerate} 

\usepackage{fullpage}  

\usepackage{setspace} 

\usepackage{hyperref} 

\usepackage{color}

\usepackage{wrapfig}

\theoremstyle{plain}

\theoremstyle{definition}

\theoremstyle{remark}

\begin{document}

\title{$90^+$ years of the Scottish Book}

\author[S. Domoradzki]{Stanis\l aw Domoradzki}
\address{Institute of History, University of Rzesz\'ow, al. Rejtana 16 c, 35-959 Rzesz\'ow, Poland}
\email{stanislawdomoradzki@gmail.com}
\author[M. Stawiska]{Ma{\l}gorzata Stawiska}
\address{AMS/Mathematical Reviews, 535 W. William St., Ste 2100, Ann Arbor, MI 48103, USA}
\email{stawiska@umich.edu}
\author[M. Zarichnyi]{Mykhailo Zarichnyi}
\address{Institute of Mathematics, University of Rzesz\'ow, 1 Profesora Stanis\l awa Pigonia Str, 	35-310  Rzesz\'ow, Poland}
\email{mzarichnyy@ur.edu.pl}

\date{\today}                                           

\maketitle

\begin{abstract}
 Inspired by the recent 90th anniversary of the Scottish Book we present some reflections about its impact.  First we discuss new areas of mathematics it helped  launch. Then we argue that it was actively used in stimulating the interests and results of junior mathematicians and students. Also, we summarize the progress  during the  decade that has passed since the publication of \cite{SBMauldin}, which contained a review of solved problems from the Scottish Book.  We also provide an overview of collections of open  problems related in one way or another to the Scottish Book. All formulations of the Scottish Book problems in English  are cited here from \cite[Chapter 6]{SBMauldin}.

\end{abstract}

\section{Introduction}

The year 2025 marked the ninetieth anniversary of the beginning of the legendary ``Scottish Book'', a notebook that was born within the walls of the ``Scottish'' coffee house in Lw\'ow and became one of the symbols of the  Lw\'ow Mathematical School. This school, formed around the figures of Stefan Banach, Hugo Steinhaus, Stanis\l aw Ulam and their colleagues, was one of the brightest centers of mathematics in the interwar Europe. Its strength lay not only in brilliant individual discoveries, but also in the joint formulation of open problems, challenges that became signposts for future generations.

The problems recorded in the ``Scottish Book'' covered functional analysis, topology, measure theory and other fields. Many of them remained open for decades, stimulating the development of new methods and directions. Some problems were solved in the second half of the 20th century, while others still inspire modern researchers. It is this practice — the collective formulation of problems — that has shown that progress in science often begins with a clearly posed question.

Today, over ninety years later, the ``Scottish Book'' is perceived not only as a historical document, but also as a living symbol: proof that mathematics is born in dialogue, in a joint search, in the ability to see a problem as the beginning of a path. The anniversary of this book is an opportunity to once again remember the Lw\'ow Mathematical School, to assess its achievements and significance for the development of modern mathematics.

An extensive historical literature is devoted to the  Lw\'ow Mathematical School. First of all, we note the book by Roman Duda \cite{Du14}, also translated into English. Of course, there is a place for the Scottish Caf\'e and the Scottish Book in it as well as in other publications. Our current article, without claiming to be complete, can serve as a supplement to existing research. In addition, 11 years have passed since the second edition of Mauldin's book \cite{SBMauldin} on the  Scottish Book, which, in particular, gave the most complete overview of the problems solved at that time or noted progress in work on unsolved problems. During these  years, mathematicians have solved a number of other problems, and this   motivated us to address these new development  while writing  this article.

\section{New and old mathematics in problems from ``The Scottish Book''}

We commonly refer to ``the" Scottish Book, but it should be noted that there are several somewhat different versions of the Scottish Book in circulation, each one with some kind of endorsement.  On the site  \cite{SBonline}, the following are accessible: a copy of the manuscript from the private collection of the Banach family (containing 193 problems, the last one unnumbered); a copy of the typescript from the private collection of the Alexiewicz family, with comments and underlining by Andrzej Alexiewicz (155 problems); a copy of the  typescript in English, translated by Stanis\l aw Ulam, from the collection of Archives of the  Library of the Mathematics Department of Wroc\l aw University (193 problems). Ulam's translation  gave rise to the publication of the Scottish Book, first as the Los Alamos Laboratory preprints (in 1957 and 1977), then as a volume edited by R. Dan Mauldin (two editions, in 1979 and 2015, the latter is \cite{SBMauldin}),  containing problems, their solutions (whenever available), extensive commentaries and other material, such as personal recollections by Ulam and Mark Kac, as well as  lectures from the 1979 conference in Denton, Texas, devoted to the Scottish Book. The second edition had also  a kind of companion conference: American Mathematical Society Special Session on the Scottish Book at the Joint Mathematics Meeting in San Antonio, Texas, on January 13, 2015 (an account can be found in \cite{Sta15}).

The readers who associate mathematics  foremost with precision may be disappointed with the varying answers to two questions: how many problems are there in the Scottish Book and what mathematical disciplines they represent? But the Scottish Book is not a mathematical object itself. It is a product of human activity and therefore susceptible to imprecision. First, the number of problems would vary because some problems consist of several parts, which may be counted  separately, especially when one part has been solved and other parts are still awaiting solution.  Second, classifying the mathematical themes present in the Scottish Book would give different results depending on the classification scheme used and on interpretation of the boundaries between different areas (which are fluid).   From 1868 the scheme of {\it Jahrbuch \"uber die Fortschritte der Mathematik} was used in classifying mathematical problems and publications. The emergence of new disciplines and growing specialization, as well as other circumstances (see \cite{Tes15} and the references therein), led to development  of a very different scheme, the Mathematics Subject Classification (MSC). This scheme was started in 1940 and was revised multiple times. 
Bearing this in mind we propose the following thematic classification of the problems in the Scottish Book (following the current MSC2020):

Set theory- 7 problems; 

combinatorics-- 5; 

group theory-- 5; 

topological groups-- 9; 

functions of a real variable-- 27; 

measure theory-- 15; 

potential theory--2; 

sequences, series, summability- 13; 

harmonic analysis-- 2; 

functional analysis-- 19; 

operator theory--4; 

convex geometry-- 9; 

general topology--25; 

topology of manifolds--16

probability--10;

others-- ...\\

Among the``others" two subsets are worth mentioning: 
Properties of polynomials--7; 
Algebraic, differential and/or elementary geometry--11. \\

Was there any major area of mathematics not represented in the Scottish Book? Let us start with  an observation of Mark Kac: 
``Returning to the Scottish Book, I would like to point out that although the
problems in it range over most of the principal branches of Mathematics, one branch
is conspicuously absent, and that is Number Theory. The reason is simple, and it is that Number Theory was not in vogue in Poland at the time. Sierpi\'nski in his
younger years (and also toward the end of his life) did important and interesting
things in Number Theory, and the Warsaw school did produce two `mutants':
A. Walfisz (who left Poland for the Soviet Union and was Professor at Tbilisi) and
S. Lubelski. There was even a serious journal, Acta Arithmetica (which continues to
this day), devoted to Number Theory, but this beautiful and important area was far
from the forefront of mathematical preoccupation in Poland before World War II." (\cite{Kac81})

The full story is more complex. We do not claim here that Kac, a participant in the Scottish Caf\'e sessions and later an author or co-author of works in probabilistic number theory (among others), made a mistake or forgot something important. Rather, we would like to point out that there may be  
further complications in thematic classification of the specific problems. While a problem may be formulated in terms belonging to some area of mathematics, its solution may use facts and techniques from other areas, or the problem itself may be closely related to problems in other areas. An interesting case is \\
  
\textbf{PROBLEM 59: RUZIEWICZ} 
 Can one decompose a square into a finite number of squares all different?\\

This looks just like yet another recreational mathematical puzzle. In fact, Ruziewicz here stated in writing a problem that he considered some years earlier in his seminar (according to \cite{SBMauldin}, he attributed it to ``mathematicians from Krak\'ow"). An analogous problem for rectangles, which  was also considered,  got solved by Zbigniew Moro\'n (1904-1971). He was a junior assistant to Ruziewicz at UJK, and after World War II he taught mathematics at schools in Wroc\l aw. He proved that the rectangle $33\times 32$ can be decomposed into  9 squares and the rectangle  $65 \times 47$ can be decomposed into  10 squares;  the solution was published in \cite{Mor25}. He had more results on the topic (including relations between decompositions of a rectangle and of a square), which he did not publish, but he  gave talks on some of them in Ruziewicz's seminar in 1925--28. He also wrote an expository paper in 1955. Different decompositions of a square were found in late 1930s. No explicit mention was made of the Scottish Book problem, although Roland Sprague's solution from 1939 (the first published one) closely resembles Moro\'n's results. In another paper (\cite{Mor56}) Moro\'n formulated  several problems in elementary number theory based on the decomposition of a square (and solved some of them). More generally, according to \cite{And20}, the problems of the decomposition of a square are related to relatively prime numbers, representations of a number as a sum of squares, the Fibonacci sequences and other number-theoretic topics. Finally, the commentary by P. J. Federico to Problem 59 in \cite{SBMauldin} discusses connections of this problem and its solutions to graph theory. So the two major areas---number theory and graph theory--- appear unexpectedly in the solutions of a problem (and its relatives) which was originally formulated without reference to them. \\

 Another problem with number-theoretic (and also probabilistic) content is the following:

\textbf{PROBLEM 152: STEINHAUS} 
A disk of radius $1$ covers at least two points with integer coordinates $(x, y)$ and
at most $5$. If we translate this disc through vectors $nw$ $(n = 1,2,3, . . .)$, where $w$
has both coordinates irrational and their ratio is irrational, then the numbers $2, 3, 4$
repeat infinitely many times. What is the frequency of these events for $n\to \infty$? Does
it exist?\\

The commentary by Jan Mycielski in \cite{SBMauldin} explains the solution using facts from diophantine approximation.\\

Graph theory was combined with probability in the following problem: 

 \textbf{PROBLEM 38: ULAM}  Let there be given $N$ elements (persons). To each element we attach $k$ others
among the given $N$ at random (these are friends of a given person). What is the
probability  $P_{kN}$ that from every element one can get to every other element through
a chain of mutual friends? (The relation of friendship is not necessarily symmetric!)
Find $\lim_{n=\infty} P_{kN}$ (=$0$ lub $1$?)."  \\
 As Ron Graham observed in his commentary  to the second edition of  \cite{SBMauldin}, this problem concerns connectivity of certain random graphs.     Hence Ulam can be regarded as the forerunner of  random graphs, over 20 years before they were formally defined  by P. Erd\H os and  A. R\'enyi (1959).\footnote{The commentary to this problem in the first edition of \cite{SBMauldin} was written by Erd\H os  and contains his solution of the problem in the case where the relation of friendship is  symmetric (the graph is undirected) and  $k \geq 2$. In the second edition,  after Erd\H os's and Graham's comments there is a detailed general solution of the problem by D. S. Jungreis. Graham does not mention Edgar Gilbert, who introduced another kind of a random graph in   1959.}  In this way Problem 38 started  a brand new area of mathematics: theory of random graphs; and, indirectly, mathematical theory of percolation, which deals with connectivity of random graphs.\footnote{In 2010 Stanislav Smirnov received the Fields medal (one of the highest honors in  mathematics)  in part for his work on percolation.}\\

As an aside, let us remark that graphs and related notions serve to describe various phenomena pertaining to science as done by human beings. In information  science and bibliometrics citation graphs are used, and collaboration graphs appear in sociology of science. The popular ``Erd\H os number" is  a distance (between mathematicians) on the Erd\H os collaboration graph.\\

The problems either originated in discussions and were formulated on the spot, or were considered by their proponents for some time beforehand.\footnote{According to Ulam, `` [The problems] 
were discussed a great deal among the persons whose names are
included in the text, and then gradually inscribed into the `book'
in ink. Most of the questions proposed were supposed to have had
considerable attention devoted to them before an `official' inclusion into the `book' was considered. As the reader will see, this
general rule could not guarantee against an occasional question
to which the answer was quite simple or even trivial." (preface to the 1957 English translation of the manuscript of the Scottish Book, \cite{SBonline} or pages xi-xiii in \cite{SBMauldin})} In the latter cases the solutions, or at least partial solutions, might have been already known. 
One such  problem is Problem 59 discussed above, another one is\\

\textbf{PROBLEM 118 FR\'ECHET}
Let $\Delta(n)$ be the greatest of the absolute values of determinants of order $n$ whose
terms are equal to$\pm 1$.  Does there exist a simple analytic expression of $\Delta(n)$ 
as a function of $n$; or, more simply, determine an analytic asymptotic expression
for $\Delta(n)$.\\
Henryk Minc in his commentary in \cite{SBMauldin} to Problem 118 discusses history of this problem. He recalls that Jacques Hadamard proved in 1893 that an arbitrary $n \times n$ complex matrix  $A$ with entries whose absolute values do not exceed   $1$ satisfies the inequality $|{\rm det} A| \leq n^{n/2}$. \\
For matrices with real entries the equality holds if and only if   $AA^T =nI_n$. Matrices satisfying this equality are called  Hadamard matrices. They found applications in coding theory. Hadamard matrices exist for $n=1,2$ or $n$ divisible by $4$. However, it is still not known if they exist for all such $n$.    Some construction go back to  J. J. Sylvester in  1867.  In 1933 Raymond Paley gave a method of constructing such matrices using finite fields for some $n < 92$. The smallest order for which the existence of a Hadamard matrix   is not currently known is $n =  668$. \\

Another new area of mathematics initiated by a problem in the Scottish Book is the theory of infinite games.\footnote{Some mathematical treatment of games can be found already in 17th century. However,  general game theory as a separate area of mathematics is not much older than the Scottish Book.  The  foundational  paper by John von Neumann was published in 1928 (\cite{Neu28}).}  Such games were considered by Stanis\l aw Mazur in \\

 \textbf{PROBLEM 43: MAZUR}\\
Prize: One bottle of wine, S. Mazur.\\
Definition of a certain game. Given is a set $E$ of real numbers. A game between two players A and B is defined as follows: A selects an arbitrary interval $d_1$, B then selects an arbitrary segment (interval) contained in $d_1$; then A in his turn selects an arbitrary segment $d_3$ contained in $d_2$ and so on. A wins if the intersection $d_1,d_2,..., \ d_n...$ [sic] contains a point of the set $E$; otherwise, he loses. If $E$ is a complement of a set of first category, there exists a method through which A can win; if $E$ is a set of first category, there exists a method through which B will win.\\
\textbf{Problem:} It is [sic] true  that there exists a method of winning for the player A only for those sets $E$ whose complement is, in a certain interval, of first category; similarly, does a method of win [sic] exist for B if $E$ is a set of first category?"\\

Problem 43 is followed by two modifications of Mazur's game, one by Ulam and another by Banach (unnumbered).  Two more variants were proposed by Banach as Problem 67. Ulam's modification of Problem 43 (known as the Ulam game, or the binary game) inspired Hugo Steinhaus and Jan Mycielski to introduce in 1962 a new axiom for set theory, called the Axiom of Determinacy, which says that the Ulam game is determined (that is, one of the two players has a winning strategy). This axiom is incompatible with the Axiom of Choice. For more details on infinite games and their influence on set theory and topology, see \cite{Tel87}.




 Theory of functions of a real variable, while modern, was already well established before the formation of the Lw\'ow School. There are 27 problems in this area in the Scottish Book,  a substantial number.  
Why are there so many problems in the area that is not considered to be a major interest of the Lw\'ow School of Mathematics, unlike functional analysis, operator theory or probability?  An explanation may be given (somewhat speculatively) that Lw\'ow mathematicians tried to get to know (confirm/refute) properties of particular functions before considering properties of whole spaces. Indeed, mathematicians often employ ``toy models",  investigating concrete but somehow typical examples to better understand, explain and foresee mathematical phenomena.  
There is also another explanation (complementary  rather than contradictory). Regarding Ruziewicz, let us recall that he got his PhD  in 1918 for a thesis that concerned functions of a real variable. His supervisor was Wac\l aw Sierpi\'nski. The areas of topology, set theory and functions of a real variable were developed in Lw\'ow even before World War I thanks to Sierpi\'nski together with Zygmunt Janiszewski and Stefan Mazurkiewicz. The three mathematicians later became the pillars of Warsaw School of Mathematics, more closely associated with the mentioned areas. Anyway,   the theory of functions of a real variable was pursued in Lw\'ow in its own right even later and   functional analysis was used to treat problems in this area. See for example Steinhaus's paper \cite{Ste29} in the first issue of  \textit{Studia Mathematica}. Besides Ruziewicz, authors of problems in the theory of functions of a real variable were J. Schauder, S. Mazur, H. Auerbach, S. Banach, Z. \L omnicki, S. Ulam, J. Schreier, H. Steinhaus,  M.Kac (affiliated with the academic institutions of Lw\'ow) , L. Infeld, A. Zygmund, W. Sierpi\'nski, A. J. Ward, S. Stoilow, E. Szpilrajn, , S. Saks and  S. Sobolev (visitors).  For some  examples (Problems 17.1 and 109), see the next section. \\

\subsection{The Scottish Book and formation of junior mathematicians}

 The youngest regular participants of the sessions in the Scottish Caf\'e were Stanis\l aw (Stan) Ulam and Marek (Mark) Kac.    Ulam got his PhD under Kazimierz Kuratowski at the Lw\'ow Polytechnics in 1933; Kac got his  PhD  under Steinhaus at the Jan Kazimierz (Lw\'ow) University in 1937. Both ultimately left for the USA on scholarships (Ulam in 1935, Kac in 1938), but visited Lw\'ow as long as they could. They contributed actively to the Scottish Book as problem-posers.  
Kac formulated Problems 126, 161, 177, 178; Ulam formulated Problems 13, 16, 17, 17.1, 18, 19, 20, 21, 29, 30, 31, 32, 33, 34, 35, 36, 37, 38, 42, a modification of 43, 68, 70, 71, 77, 99, 101, 102, 110, 115, 139, 140, 141, 142, 145, 146, 164, 165, 166, 167, 168 (alone),  2, 3, 11, 40, 100 (with Banach), 10.1 (with Mazur, Auerbach and Banach), 22, 95, 98, 103, 116, 163 (with Schreier), 62, 63, 69, 109, 144 (with Mazur), 90, 114 (with Auerbach), 94 (with Z. \L omnicki), 96 (with Kuratowski), 108 (with Banach and Mazur)
and solved Problem 9.\\
Zbigniew \L omnicki (1904-1994) did not attend regularly, since he worked as an actuary after finishing his studies in mathematics and physics at the Lw\'ow University and in 1933 moved to Warsaw for a job at Zak\l ad Ubezpiecze\'n  Spo\l ecznych [Social Insurance Institution]. But he co-authored Problem 94 with his mathematical collaborator Ulam. We should also mention here  J\'ozef Marcinkiewicz (1910-1940) (\cite{DPB}, \cite{Ma11}, \cite{MJ11}). Educated at the University of Wilno (Vilnius), he got his master's degree in 1933 and  PhD in 1934. Then  he spent the academic year 1935/36 at the Lw\'ow University, supported by a scholarship of the National Cultural  Fund [Narodowy Fundusz Kultury]. Despite his young age, he was already an established mathematician. In Lw\'ow he was an assistant to Banach   and discussed mathematics with others, in particular with Juliusz Schauder and  Stefan Kaczmarz (he published a joint paper with the latter).  Marcinkiewicz  attended the Scottish Caf\'e and  solved several problems from the Scottish Book: 83, 106 (his solution needed interpretation), 131, and formulated Problem 124. Upon his return to Wilno he tried to establish a similar tradition of mathematical discussion among professors, assistants and senior students in a caf\'e setting, but it did not gain ground. According to Leon Je\'smanowicz, then an assistant in Wilno (\cite{Je11}), the reasons for failure were high prices, early closing time and  finally the temporary closure of the university because of student riots.\\

Some other junior mathematicians affiliated at some point with Lw\'ow, but not  necessarily trained within the Lw\'ow school of mathematics, were inspired by the problems entered in the Scottish Book and worked on their solutions with success.  The first example is Andrzej Turowicz (1904-1989). He graduated from the  Jagiellonian University (Krak\'ow) with master's degree in mathematics in 1928, then he spent several years teaching mathematics in  high schools  in Krak\'ow and Mielec. In 1937, at the initiative of Antoni \L omnicki, he was appointed  a senior assistant in the Chair of mathematics of Lw\'ow Polytechnics. He was active in the mathematical life of Lw\'ow and collaborated with mathematicians from the Lw\'ow School of Mathematics.  Together with Stefan Kaczmarz (1895-1939) he solved the following problem posed by Mazur and Ulam:

\textbf{PROBLEM 109: MAZUR, ULAM}

October 16, 1935

Given are $n$ functions of a real variable: $f_1, . . . , f_n)$. Denote by $R(f_1, . . . , f_n)$ the
set of all functions obtained from the given functions through rational operations
(expressions of the form
\[
\frac 
{\sum a_{k_1...k_n} f_1^{k_1}\cdots f_n^{k_n}}{\sum b_{k_1...k_n} f_1^{k_1}\cdots f_n^{k_n}} \biggr). 
\]
Must there always exist, in the set $R$, a function $f$ such that its indefinite integral
does not belong to the set $R$?\\
An analogous question in the case where we include in the set $R$ all the functions
obtained by composing functions belonging to $R$.\\

Their solution was published in the leading journal \textit{Studia Mathematica}, based in Lw\'ow, as the paper \cite{KaTu39}. The addendum in \cite{SBMauldin}, following  Ulam's typewritten translation, says that ``An affirmative answer for the first question was found by Docents,
Dr. S. Kaczmarz and Dr. A. Turowicz."\footnote{In the Polish copy of the manuscript of the Scottish book as well as in the Polish typescript from the Alexiewicz family collecion, available at \cite{SBonline}, there is no addendum to this problem.} It is dated March 1938, but must have been added much later. Turowicz obtained his PhD degree only in 1946, at the Jagiellonian University the on the basis of the thesis ``On continuous and multiplicative functionals" , under the supervision of  Tadeusz Wa\.zewski. \footnote{For more about the fascinating life and achievements of Turowicz,  a Benedictine monk after World War II, see \cite{CP88} and \cite{DoSt15}.} 
Turowicz's PhD  thesis (published as \cite{Tu48})  answered a question posed by Stefan Banach and Meier Eidelheit. Even though both Banach and Eidelheit entered many problems in the Scottish Book, this was not one of them.  Here is how Turowicz recalled the origin of  his thesis work:

\textit{I had an incident with Banach like this:  for  a meeting of the mathematical society, I proposed a talk on multiplicative and continuous functionals (my proposal was in spring 1939; I finished the work after the war). I am delivering the talk and Banach enters the room, slightly late, with an incredibly sullen face. I noticed that  Banach was angry. He listened with extreme attention and his face changed. When I finished, Banach took the floor and said: ``I also dealt with this problem; you did it in a totally different way, and you did it well." I received his opinion with gladness. The next day after this meeting, Sto\.zek (who was not at the meeting) asked: ``Was Banach there?" [I said]  ``Yes, he was." [He said] ``I did not want to scare you in advance; he was very angry when he found out what you were to talk about. He said: `I am dealing with this; I must have told someone, and [now] Turowicz is presenting it as his own.'" Banach came with the intention of giving me a hard time. Luckily the idea of the proof was completely different [from his], therefore [he] praised me [and] did not make a scene. (...) Since then, Banach was very friendly towards me.}
(\cite{TuAU}, cassette 2b), translated from the Polish by M. Stawiska)
This incident supports Ulam's observation about problems being given ``considerable attention" before being shared, but we will  never know for certain whether Banach ever intended to enter the problem in question in the Scottish Book. 

Another example of a Scottish Book problem inspiring work of a junior mathematician is 

Problem 17.1 (Ulam): Let $f$  be a continuous function defined for all $0 \leq x \leq 1$. Does there exist a perfect
set of points $C$ and an analytic function  $\varphi$  so that for all points of the set $C$ we have
$f \equiv \varphi$?\\

It was solved by Zygmunt Zahorski (1914-1998).  Zahorski studied mathematics at the University of Warsaw in 1934--38, then taught at Air Force Cadet School and was preparing a PhD thesis under direction of Stefan Mazurkiewicz. When World War II broke out in 1939, he fled to Lw\'ow. In  1939-1941 he was an assistant to  Banach and an  aspirant (an equivalent of a doctoral candidate) at the Lw\'ow Uniwersity. According to \cite{Li00}, he was told about  Ulam's problem by Banach. He solved it  in 1941, at the end of the first Soviet occupation of Lw\'ow and the beginning of the Nazi German one.  The publication was of course impossible then; it happened only after the war  (\cite{Zah47}). The solution is negative\footnote{In \cite{Zah47}, Zahorski reformulates the problem, attributing it also to Ulam. The positive solution of the reformulated problem implies negative solution to Problem 17.1.}, following from a series of lemmas,  and appears along with solutions of several other problems and new related open questions. Zahorski characterized it in \cite{Zah92}, tongue in cheek, as follows: {\it ``The three-minute proof by contradiction I  am leaving to the reader, it is enough to pick a function such that $r(x)=0$ for every $x$, but the construction of such a function may take several hours, and [for] the whole condition a month is enough, for me, or it may fail."} Despite these complications--or maybe because of them?--he considered the solution to Problem 17.1 to be one of his main achievements. Also in \cite{Zah92} he confessed that he worked briefly on another Scottish Book problem, the famous Problem 153 of Mazur, ultimately solved by Per Enflo: {\it ``Sometime right after the war I tried to solve Mazur's problem--I found it difficult and it seemed artificial, I did not know that it is so closely connected to a  natural problem of convergence and I abandoned my attempts---otherwise I would still try (difficulties just attract me)."}\footnote{For more on the life  and achievements of Zahorski, see \cite{DuTe}, \cite{Zah86}.}\\

 From the times of the first Soviet occupation of Lw\'ow we should mention also the case of Mark Vishik (1921-2012) and Vladyslav Lyantse (1920-2007), who  then were  students of mathematics at the Lw\'ow University. Conferences specially for students started to be organized and in April 1941 Vishik gave a talk at such a conference, supervised by senior mathematicians such as Banach and Edward Szpilrajn (a refugee from German-occupied Warsaw, later  known as Marczewski).\footnote{For more information on this conference and its participants, see \cite{MP}.}  It concerned his joint work with Lyantse, partially solving Problem 192 form the Scottish Book, posed by Bronis\l aw Knaster (another refugee from Warsaw) and Szpilrajn (\cite{Kal}, \cite{SCKMS}). They did not publish this solution (perhaps could not, because of the start of the German-Soviet war soon after), but it was mentioned  in the remarks in the original Scottish Book. Szpilrajn also mentioned it in his own paper published right after the war \cite{Sz45}).\\

\textbf{PROBLEM 192: B. KNASTER, E. SZPILRAJN}
May, 1941\\
Definition. A topological [sic] $T$ has the property $(S)$ (of Suslin) if every family of
disjoint sets, open in $T$, is at most countable.
Definition. A space $T$ has property $(K)$ (of Knaster) if every noncountable
family of sets, open in $T$, contains a noncountable subfamily of sets which have
elements common to each other. (...)

Problem (B. Knaster and E. Szpilrajn). Does there exist a topological space (in
the sense of Hausdorff, or, in a weaker sense, e.g., spaces of Kolmogoroff) with the
property $(S)$ and not satisfying the property $(K)$?
(...)
Problem (E. Szpilrajn). Is the property $(S)$ an invariant of the operation of
Cartesian product of two factors?\\
Remarks: (...)
(5) E. Szpilrajn proved in May 1941 that the property $(K)$ is an invariant of the
Cartesian product for any number of factors and B. [sic] Lance and M. Wiszik
verified that if one space possesses property $(S)$, and another space has property
$(K)$ then their Cartesian product also has property $(S)$.\\

In \cite{Ul}, Ulam views the caf\'e sessions as extensions of the meetings of the Lw\'ow branch of the Polish Mathematical Society, even before the Scottish Caf\'e became the preferred place: \textit{``For our story, the Caf\'e Roma was, in the beginning, the more important of the two
coffee houses. It was there that the mathematicians first gathered after the weekly
meetings of our local chapter of the Polish Mathematical Society. The meetings
were usually held on Saturday in a seminar room at the University--hence close
to the Caf\'es. The time could be either afternoon or evening. The usual program
consisted of four or five ten-minute talks; half-hour talks were not very common,
and hour-long talks were mercifully rare. There was of course some discussion at
the seminar, but the really fruitful discussions took place at the Caf\'e Roma after the
meeting was officially over."} It would be interesting to compare the topics of the Polish Mathematical Society meetings with the problems entered into the Scottish Book.  This will have to wait until the minutes of these meetings are accessible (if they exist). The reports published in Annales de la Soci\'et\'e Polonaise
de Math\'ematique (the journal of the Polish Mathematical Society) for the years 1935-1938 submitted by the Lw\'ow branch are incomplete. There is no information about the meetings between July 1935 (when Problem 1 was entered by Banach) and June 1938. The reports for the first half of 1935 and those for 1938 are lacking abstracts, so  is hard to determine whether the topic of the talks  are related to any  problems entered in the Scottish Book by the speakers. For the second half of  1938, only 4 meetings were reported, with talks by Jean  Leray, S. Ulam, S. Mazur and M. Kac. Kac entered 2 problems (177 and 178) in September 1938, before his October talk.\\


\section{Rewards of doing mathematics}

For solving some problems form the Scottish Book, their authors awarded prizes, mostly drinks and food to be consumed right away. A notable exception was the live goose proposed by Stanis\l aw Mazur, which eventually, after being presented in 1972 to Per Enflo for the solution of Problem 153, also became food (\cite{Ku}). Sometimes different prizes were offered according to whether the problem would be solved  positively or negatively.  In connection with his Problem 152 (November 6, 1936), Hugo Steinhaus offered the following prizes: 
``For the computation of the frequency: 100 grammes of
caviar. For a proof of the existence of frequency, a small beer. For counter example [sic]:
A demitasse."

However, rewards for solving scientific problems were not an  invention of the Lw\'ow School of Mathematics, but predate it substantially. Many institutions (academies of sciences, learned societies, universities and others) organized competitions and offered prizes to determine the ``best" work, sometimes on a single (yet broad enough) topics. Among the many such competitions let us highlight the 1915-1918 Grand Prix des Sciences Math\'ematiques by the French Academy of Sciences in Paris on the topic of iteration of functions of  complex variables (one or several). The papers of Gaston Julia, Samuel Latt\`es and Salvatore Pincherle entered into this competition (along with the contemporaneous  works of Pierre Fatou, who did not compete) are cornerstones of the modern theory of holomorphic dynamics (\cite{Au}, \cite{Rosa}). Ultimately Julia won the 3000 francs.\\

Prize competitions could also focus on individual problems. 
G\"osta Mittag-Leffler (founder and editor of Acta Mathematica) initiated an international mathematics competition in 1885 to celebrate the 60th birthday of King Oscar II. Four unsolved problems were proposed by the commission (Karl Weierstrass, Charles Hermite). The prize was gold medal + 2,500 Swedish crowns + publication in Acta Mathematica. After winning the competition in 1889, Henri Poincar\'e's results were prepared for publication. When a serious error in his argument for stability was discovered during printing, Poincar\'e immediately admitted the error and returned the prize money to the organizers, covering not only the prize but also the printing costs already incurred.\footnote{See \cite{BG} for an exhaustive analysis.}\\

Paul Friedrich Wolfskehl (1856-1906), a physician interested in mathematics, designated in his will  100,000 marks as the prize for the proof of Fermat's Last Theorem (to be awarded by  the Royal Society of Science in
G\"ottingen). Andrew Wiles received the prize in 1997,  reduced by inflation to about 75 000 DM (\cite{Bar}).\\

 The tradition of offering  small prizes for  solutions of  mathematical problems was  continued after World War II. Mark Kac in 1981 offered ten Martini cocktails for the solution of the problem  whether  the spectrum of the almost Mathieu operator is a Cantor set. The problem  was solved positively in 2009 by A. Avila and S. Jitomirskaya, with the solution published in one of the most prestigious mathematical journals (\cite{AJ}). 

It is also natural here to mention  Paul Erd\H{o}s, who was known for his own  tradition of awarding monetary rewards for solving mathematical problems he proposed. He assessed the difficulty of the problem on his own scale. The amounts ranged from \$25 (for relatively simple problems) to \$10,000 (for particularly difficult problems).\\

In contrast to this, a million dollars is offered for each of the well-known Millennium problems chosen by the Clay Institute. The same amount is offered by the American banker Andrew Beal for proving or disproving his conjecture: If
$ A^{x}+B^{y}=C^{z}$,
where $A, B, C, x, y$, and $z$ are positive integers with $x, y, z > 2$, then $A, B$, and $C$ have a common prime factor.\\

The Lw\'ow prizes did not generally cost much money. Perhaps in the time of deprivation during the first Soviet occupation\footnote{Here is how Steinhaus described the difficulties with provisions: ``These same authorities
kept promising that \textit{vsyo budet}, but the shortages persisted: among much else,
soap, sugar, kerosene, and fuel for heating were in extremely short supply—and it
was a harsh winter."\cite{Ste15}}  the prize for solving Problem 184 (posed in 1940 by Stanis\l aw Saks, yet another refugee from Warsaw) of a kilogram  of fatback\footnote{In Polish ``s\l onina". Mistranslated in \cite{SBMauldin} and earlier editions as ``bacon".} would be particularly welcome, even if not so cheap.  Many problems were not assigned a prize. In one instance a prize was offered by someone other than the problem's author: Samuel Eilenberg offered a prize of one bottle of wine for the solution of part (a) of Problem 77 posed by Ulam.\\

However, the prizes themselves were not the ultimate  goals. This is a broader phenomenon. For example, 
Alexandre Grothendieck refused the Crafoord Prize (1988) and 
Grigori Perelman refused the Fields Medal (2006) and the Millennium Prize  (2010).


To draw some general conclusions, it is worth analyzing the phenomenon of the Scottish Book in terms the concepts developed by Polish-Israeli physician Ludwik Fleck (a collaborator of Steinhaus in Lw\'ow).

According to Fleck \cite{Fleck}, a thought collective (``Denkkollektiv'' in German)  is a community of researchers who collectively interact to produce or develop knowledge, using a shared system of cultural customs and knowledge acquisition.

Scientific knowledge production can be defined primarily as a social process that depends on previous discoveries and practices in a way that limits and conditions new ideas and concepts.

Fleck called this shared collection of pre-existing knowledge a ``denkstyle'' or style of thinking, and formulated a comparative epistemology of science using these two ideas.

By interacting socially in the process of knowledge production, researchers create shared concepts and practices that they use to discuss and debate each other's ideas and discoveries. 

The Scottish Caf\'e and the Scottish Book were, in a sense, platforms for social interaction among the participants of the Lw\'ow Mathematical School. The juniors were considered intellectual peers. The hospitality was  extended to visitors from other  centers, even to representatives of Soviet academic institutions during the first Soviet occupation of Lw\'ow. \footnote{Problems were entered into the Scottish Book by N. N. Bogolyubov (183), P. Alexandroff (187), S. Sobolev (188), A. F. Fermant (189) and L. Lusternik (190).}  The  role of prizes was largely symbolic, but they emphasized the friendly atmosphere that prevailed among  Lw\'ow  mathematicians. \\

 \section{During the last decade}
 
 The aim of this setion is to rewiev the results concerning the problems of the Scottish book and obtained during the last decade, i.e. after publication of \cite{SBMauldin}. As we will see, there is significant progress, even to the point of complete solution, of some problems. In some cases, we also provide additional information to the comments from \cite{SBMauldin}.
 
 Given the vastness of the subject matter of the Scottish Book, the material in this section is inevitably far from being complete.
 
 Set-theoretic (general) topology is related to functional analysis, which in a broad sense involves the study of topological and algebraic structures and their interplay. This is probably one explanation for the fact that a significant number of the problems formulated in the Scottish book concern topology. Problems related to the theory of topological groups and semigroups should be added here.

 PROBLEM 1: BANACH 
 
 July 17,1935 
 
 (a) When can a metric space [possibly of type (B)] be so metrized that it will become complete and compact, and so that all the sequences converging originally should also converge in the new metric? 
 
 (b) Can, for example, the space $c_0$ be so metrized? 
 
 It is remarked in \cite{BaPl} that in te modern terminology te problem sounds as follows:
 
  When does a metric (possibly Banach) space $X$ admit a condensation (i.e. a bijective continuous map) onto a compactum (= compact metric space)?
  
  Problems concerning  condensations of topological spaces belong to set-theoretic topology. It is possible that the initial problems from the Scottish Book were discussed for years at the Lw\'ow Mathematical School. The article \cite{BaPl} suggests that the roots of Problem 1 lie in Sierpinski's article \cite{Sie}. Later, the problem of compactifications of topological spaces into compact spaces was considered in the Alexandrov school. \\
  
  Mauldin's commentary notes that the solution to the problem for separable Banach spaces follows from the Anderson-Kadets theorem on the homeomorphism of such spaces to the space $\mathbb R^\infty$, and for the latter such compaction is easy to construct, starting from the  condensation of $\mathbb R$ onto a compact set.
  
  However, Mauldin does not say anything about non-separable Banach spaces. The main result of the paper \cite{BaPl} is the theorem that the Hilbert space of density $\omega_1$ condenses onto the Hilbert cube $Q$. Since all Banach spaces of the same density are homeomorphic, this provides a complete solution of the problem in the class of Banach spaces of density $\omega_1$.
  
  Also, in \cite{BaPl} a relatively simple solution of Problem 1 for Banach spaces of density $\mathfrak c$ is given. The question of densities between $\omega_1$ and continuum remained open. 
  
  In \cite{Banak} Taras Banakh announced the following results:
  
  (1) It is consistent that the continuum is arbitrarily large and every infinite-dimensional Banach space of density $\leq c$ condenses onto the Hilbert cube. 
  
  (2) It is consistent that the continuum is arbitrarily large and no Banach space of density $\aleph_1 < d(X) < c$ condenses onto a compact metric space.

   Here ``consistent'' means ``consistent with Zermelo-Fraenkel set theory''. The fact that the answer to the Banach problem depends on the axiomatics of set theory makes it similar to Cantor's Continuum Hypothesis.\\

  The ``metric'' part of Banach's problem was solved by Pytkeev \cite{Pyt} for the class of separable absolute Borel spaces. 
  
  
  
  
  
 

For related developments, see \cite{Osi}, \cite{Osi1}.\\
 
 PROBLEM 7: MAZUR, BANACH
 
 Are two convex infinite-dimensional subsets of a Banach space [of type (B)]
 always homeomorphic?
 
 In this formulation, the solution of the problem is trivial; as Taras Banakh noted \cite{Banak}, Problem 7 should have been the problem of classifying closed convex sets in Banach spaces up to homeomorphism. In this formulation, the problem turned out to be important for the development of infinite-dimensional topology.
 
 The list of references in V. Klee's commentary ends with year 1975. We are going to extend it by mentioning some of the most important results in this direction. 
 
 H. Toru\'nczyk \cite{Tor} proved that any Fr\'echet space is homeomorphic to a Hilbert space.
 
Dobrowolski and Toru\'nczyk \cite{DoTo} extended the mentioned Klee's classification: 

Each separable closed convex subset $C$ of a Fr\'echet space is homeomorphic to $[0,1]^n \times [0,1)^m \times (0,1)^k$ for some cardinals $0 \le n,k \le\omega$ and $0 \le m \le 1$. In particular, $C$ is homeomorphic to the separable Hilbert space $\ell^2$ if and only if $C$ is not locally compact.

Finally, the following  result by Banakh and Cauty \cite{BaCa} can be regarded as a complete solution to the (modified) Problem 7 in the class of Fr\'echet spaces:

Each closed convex subset $C$ of a Fr\'echet space is homeomorphic to $[0,1]^n \times [0,1)^m \times \ell_2(\kappa)$  for some cardinals $0 \le n \le\omega$, $0 \le m \le 1$ and $\kappa\ge 0$. In particular, $C$ is homeomorphic to an infinite-dimensional Hilbert space if and only if $C$ is not locally compact.\\
 
 
 
 
 
 
 PROBLEM 19: ULAM Is a solid of uniform density which will float in water in every position a sphere?

 Since the second edition of the Scottish Book, some new papers have been published in this direction. In \cite{HSW}, a new family of convex bodies related to Problem 19 was introduced and investigated.
 
 A short affirmative solution of Problem 19 in the class of symmetric convex bodies in $\mathbb R^n$ is given in \cite{FSWZ}.
 
  In 2022, the problem was solved negatively \cite{Rya22}. The solution cannot be symmetric, but it is a body of revolution.
  
  In \cite{Zav}, the author presents natural affine counterparts of the classical theorems on Problem 19 due to H. Auerbach.
 
 For discussion of related developments, see \cite{ASG}, \cite{FSWZ}, \cite{LWYZ}, \cite{Mon24}, \cite{Rya23}.

 PROBLEM 35 by Stanis\l aw Ulam:
 
 Is projective Hilbert space (that is to say, the set of all diameters of the unit sphere
 in Hilbert space metrized by the Hausdorff formula) homeomorphic to the Hilbert
 space itself?
 
W. Holsztynski commented that the answer is no because Hilbert space is simply connected and projective Hilbert
 space is not simply connected (because it has a double covering). Definitely, this immediately leads to the question of determination of topological type of projective Hilbert space.  
 
 This is closely related to Problem 36, also formulated by Ulam and which can be reformulated as follows: is there a retraction of a closed unit ball of a Hilbert space to its boundary? The answer to problem 36 is positive and, moreover, the unit sphere of a Hilbert space is homeomorphic to the entire space. It follows that a projective Hilbert space is a Hilbert manifold. Also, since the group $\mathbb Z/2\mathbb Z$ acts freely on the  unit sphere of a Hilbert space, the projective Hilbert space is a space $K(\mathbb Z/2\mathbb Z,1)$. Since Hilbert manifolds are classified by their homotopy type, the Hilbert projective space is homeomorphic to the product of an infinite telescope over an infinite projective space onto the Hilbert space. This is what the complete solution of Problem 35 should look like.

 PROBLEM 42 by Stanis\l aw Ulam:
 
 To every closed, convex set $X$, contained in a sphere $K$ in Euclidean space, there
 is assigned another convex, closed set $f(X)$, contained in $K$, in a continuous manner
 (in the sense of the Hausdorff metric); does there exist a fixed point, that is to say, a
 closed convex $X_0$ such that $f(X_0) = X_0$?
 
 Actually, this problem concerns the hyperspace of compact convex subsets of $K$. Nadler, Quinn, and Stavrokas \cite{NKS} proved that this hyperspace is homeomorphic to the Hilbert cube $Q$. The positive answer to the problem is now precisely the fixed point theorem for $Q$ (infinite-dimensional extension of Brouwer fixed point theorem).

 PROBLEM 54 by Juliusz Schauder: 
 
 A convex, closed, compact set $H$ is transformed by a continuous mapping $U(x)$
 on a part of itself. $H$ is contained in a space of type (F). Does there exist a fixed
 point of the transformation?
 
 It is known that every map of a convex compact finite-dimensional set has a fixed point (Brouwer). Thus, Schauder's problem concerns the existence of infinite-dimensional analogues of Brouwer's theorem. 
 
 In 1930, Schauder published a paper \cite{Sch} on the infinite-dimensional fixed point theorem. However, there was a gap in the proof of the main result: it only worked for Banach spaces. In 1934, Tychonoff \cite{Tych} proved a fixed point theorem for compact convex subsets in locally convex spaces.
 
Since the appearance of the problem in the Scottish book, many fixed point results have appeared. We will only mention here the Kakutani-Glicksberg fixed point theorem for multivalued mappings \cite{Glic}. It is important and popular due to its applications in equilibrium theory and game theory.

However, in its general formulation, Schauder's problem remained unsolved. In 2001, Robert Cauty published an article \cite{Cau} arguing that the problem had a positive solution. Cauty's arguments were expanded and explained by Tadeusz Dobrowolski \cite{Dobr}. However, it soon became clear that this approach did not yield a correct proof. Cauty published the paper \cite{Cau1} in which he developed a theory of algebraic absolute neighborhood retractions, which he believed provided a complete solution to Schauder's problem.

Today, it cannot be assumed that Cauty's proof is accepted by the entire mathematical community (see \cite{KST}, which cites only \cite{Cau} and \cite{Dobr}, but not \cite{Cau1}). In any case, as \cite{KST} says, ``it is desirable to come up with much simpler arguments than Cauty's''.

Some new results related to Problems 73, 74 are obtained in \cite{ChaSa}.

PROBLEM 107: STERNBACH
Does there exist a fixed point for every continuous mapping of a bounded plane
continuum E, which does not cut the plane, into part of itself? The same for
homeomorphic mappings of E into all of itself

The problem remains open today. Some partial and related results can be found in \cite{Roba,Boro, BoKu}   and \cite{BOT}.

PROBLEM 138 by Samuel Eilenberg: Any compact convex set located in a linear space of type ($B_0$) is an absolute retract. 

(Recall that a linear space of type ($B_0$) is precisely a Fr\'echet space, (i.e. a complete metrizable topological vector space.) This problem gave rise to a whole series of studies on the theory of absolute retracts and the related theory of absolute extensors. In particular, the fundamental Dugundji theorem states that every nonempty convex subset of a locally convex metrizable topological vector space is an absolute retract.

One of the most striking results is that not all linear topological spaces are absolute retracts \cite{Cau2}.

PROBLEM 151: WAVRE

November 6, 1936; 

Prize: A ``fondue” in Geneva; Original manuscript in French

Does there exist a harmonic function defined in a region which contains a cube
in its interior, which vanishes on all the edges of the cube? One does not consider
$f \equiv 0$.

The affirmative answer is given in \cite{Vas25}.

PROBLEM 157: WARD

March 23, 1937; 

Prize: Lunch at the ``Dorothy''

The anonymous comment to this problem asserts that the problem was already solved in 1975 by  Richard O’Malley. In \cite{BHR}, it is remarked that the affirmative solution to Problem 175 is an immediate consequence of a result of Daniel Ornstein \cite{Orn}.

PROBLEM 155: MAZUR, STERNBACH 

November 18, 1936 

 Given are two spaces $X, Y$ of type (B), $y = U(x)$ is a one-to-one mapping of the space $X$ onto the whole space $Y$ with the following property: For every $x_0 \in X$ there exists an $\epsilon > 0$ such that the mapping $y = U(x)$, considered for $x$ belonging to the sphere with the center $x_0$ and radius $\epsilon$, is an isometric mapping. Is the mapping $y =U(x)$ an isometric transformation? This theorem is true if $U^{-1}$ is continuous. This is the case, in particular, when $Y$ has a finite number of dimensions or else the following property: If $\|y_1 +y_2\| = \|y_1\| + \|y_2\|$, $y_1\neq 0$, then $y_2 =\lambda y_1$, $\lambda \ge 0$.
 
 In modern terminology, this problem can be reformulated as follows:
 
 Let $U$ be a bijective map between Banach spaces. Suppose $U$ is locally isometric. Is it true that $U$ is an isometry?
 
No comment to this problem is given in \cite{SBMauldin}. For the separable Banach spaces, the problem is solved in the affirmative in \cite{Mori}. Another partial solution is obtained in \cite{Basso}.

Let $X$ and $Y$ be Banach spaces and $U: X \to Y$ a local isometry. Then $U$ is an isometry. Here, a map f$: X \to Y$ is said to be a local isometry if every $x \in X$ has an open neighborhood $V$ such that $f|V$ is an isometry onto an open subset of $Y$.



PROBLEM 175 by Karol Borsuk.

(a) Is the product (Cartesian) of the Hilbert cube $Q$ with the curve which is shaped like the letter $T$, homeomorphic with $Q$? (b) Is the product space of an infinite sequence of letters $T$ homeomorphic to $Q$?

Recall that the Hilbert cube $$Q=I\times I\times I\times \dots $$ ($I$ stands for the unit segment)
is the natural counterpart of the $n$-dimensional cube $I\times I\times I\times\dots\times I$ ($n$ factors). 

Informally speaking, Borsuk's problem was asking whether an infinite product of segments could ``absorb'' a singular point in the letter T. In H. Torunczyk's extensive commentary on this problem (which was positively solved by Anderson), its great importance for the development of the geometry of the Hilbert cube was noted. In fact, one could say more: the topology of infinite-dimensional manifolds grew from several sources, one of which is Borsuk's problem.

Some of the problems from the Scottish Book are associated with the name of J\'ozef Schreier (problems 4, 11, 22, 29, 48, 67, 68, 95, 96, 98, 103, 111, 116). In the article \cite{GuP} the authors give an overview of the results related to Schreier's doctoral thesis, as well as some problems from the Scottish Book. This applies in particular to Problem 96 formulated by Ulam:

Can the group $S_\infty$ of all permutations of integers be so metrized that the group
operation (composition of permutations) is a continuous function and the set  $S_\infty$ 
becomes, under this metric, a compact space? (locally compact?)

There was Addendum of November, 1935 by Schreier and  Ulam: One cannot metrize this group in a compact way. However, the locally compact case remained unsolved till 2011, when it was resolved negatively in \cite{BaG}.

\textbf{PROBLEM 184 S. SAKS} 
Prize: one kilo of bacon.
A subharmonic function $\phi$ has everywhere partial derivatives $\partial^2\phi/\partial x^2$, $\partial^2\phi/\partial y^2$. Is it true that $\Delta\phi\ge 0$? 

Remark: it is obvious immediately that $\Delta\phi\ge 0$ at all points of continuity of $\partial^2\phi/\partial x^2$, $\partial^2\phi/\partial y^2$, therefore on an everywhere dense set.

\medskip

The following result is proved in \cite{Ere}:

Let $u$ be a subharmonic function of two variables whose first partial derivatives exist on the coordinate axes and $u_{xx}$, $u_{yy}$ exist at the origin. Then $u_{xx}(0,0)+u_{yy}(0,0)\ge0$.

This gives the affirmative solution to this problem.


\section{Continuing traditions}

Several other collections of open mathematical problems are related to the Scottish Book in various ways (historically, geographically, or ideologically). We will mention some of them in this section.

\subsection{``The New Scottish Book''}

The most important sequel to the Scottish Book is called the New Scottish Book. It was created by mathematicians from Wroc\l aw. In fact, a brief history of the New Scottish Book constitutes the content of  \cite[Chapter 8]{SBMauldin}, and \cite[Chapter 9]{SBMauldin} contains selected problems from the book. It is very symbolic that the first problem was entered by Hugo Steinhaus, who also  entered the last problem in the original Scottish Book.

The authors do not provide exact statistics, but at first glance it is clear that most of the problems belong to topology or related areas. Some problems are placed with answers or remarks. Prizes are awarded for solving some problems. Some of prizes are intangible, such as ``an opportunity to tell two jokes of
moderate length'' for solving some problems by B. Knaster.

\subsection{``The Lviv Scottish Book''}

The premises where the Scottish Book was located changed its purpose several times after the war. For a while, there  was even a bank there. In 2014, there was a restaurant and bar there again. The owners brought back the old name Szkocka. As an attraction, the restaurant kept a copy of the Scottish Book and soon visitors were writing down on the blank pages of the copy their impressions of visiting the restaurant, but also -- unexpectedly -- new mathematical problems. Soon special notebooks were allocated for the problems and so the new Scottish Book appeared. A detailed history, as well as a list of some of the problems, can be found in \cite{BBHPZ}.

``The Lviv Scottish book'' is a collective participant in the collective blog and online community of mathematicians mathoverflow. Among the community members are many famous mathematicians from different countries of the world. 

One of the problems formulated by Polish graduate students who visited Lviv was solved by the American mathematician Terence Tao, a Fields Medal winner. Since, in full accordance with the traditions of the Lviv Mathematical School, the authors of the problem announced a prize for its solution -- a bottle of drinking honey, it was decided to present it during Terence Tao's arrival in Warsaw. At the same time, Terence Tao received a special award from the modern Lviv Mathematical Society, a sculpture of a goose.

\subsection{``The Nonarchimedean Scottish Book''}

Recall that a metric $d$ is called non-Archimedean if it satisfies the strong triangle inequality $d(x,y)\leq \max\{d(x,z),d(z,y)\}$. The $p$-adic metric is an example of a non-Archimedean metric. About ten years ago, Kiran S. (Sridhara) Kedlaya from the UC San Diego initiated the so-called ``The Nonarchimedean Scottish Book'' (see \cite{Ked}). 

\subsection{``Open Problems in Topology''}

A significant part of the problems in the Scottish Book are related to set-theoretic topology. In this sense, \cite{OPT1,OPT2} can be considered a continuation of the Scottish Book. We should also note a certain geographical connection: two chapters in the book were written by modern Lviv mathematicians. 

One of the features of the ``Open Problems in Topology'' project was the periodic publication of solved problems in the journal Topology and its Applications. As far as we know, such publications have been discontinued  some time ago.

\subsection{``Ksi\c a\.zeczka Szkocka" [The Scottish Booklet]} In the 1970s and 1980s students of mathematics at Jagiellonian University (Krak\'ow, Poland) associated with Ko\l o Ma\-te\-ma\-ty\-k\'ow Student\'ow UJ [Circle of Mathematicians Students of UJ]  maintained a notebook where mathematical problems were entered, some with prizes. The following account can be found on the Circle's website: 

\textit{``In the year 1978 the Circle, following the example of the masters from Lw\'ow, opened the ``Scottish Teahouse". There students were posing problems and setting prizes for solutions. Offered were a bar of chocolate, a Pepsi, a lantern battery; today it all sounds banal, but then these  prizes were practically unattainable in stores. The problems were recorded in a special ``Scottish Booklet".  The greatest fame was achieved by problem no. 48, later called the isometry problem, posed in the year 1979 by Edward Kania. For a year students and assistants struggled with the problem; in October 1980 the solution was given by the beginning student S\l awomir Ko\l odziej."}(\cite{Kos}, translated from the Polish by M. Stawiska).


 It was a strictly local initiative, which does not seem to be continued. It was not attached to an actual teahouse. The name, purposely related to the great Lw\'ow predecessor,  was given to a series of informal meetings at which tea was drunk (coffee being harder to obtain). Again,  the Scottish Caf\'e was treated as a symbol. \\

To conclude, let us  cite the journalist J\'ozef Mayen \cite{May34}, who first described the intellectual communities associated with the Scottish Caf\'e: \textit{``The coffeehouse is a meeting which one does not schedule, but which one  comes to! \textit{Scilicet:} A meeting not only with others, but also a unique opportunity for free meeting with oneself- a unique opportunity of unconstrained ``coming to oneself".}


\section{Acknowledgements} The authors thank Maciej Klimek for bringing the reference \cite{Kos} to their attention.

\end{document}